\documentstyle{article}
\begin{document}
\begin{center}
{\large\bf SOLUTIONS OF FRACTIONAL REACTION-DIFFUSION EQUATIONS\\ 
IN TERMS OF MITTAG-LEFFLER FUNCTIONS}\\[0.5cm]
{\bf R.K. Saxena}\\
Department of Mathematics and Statistics, Jai Narain Vyas University\\
Jodhpur 342004, India\\[0.5cm]	
{\bf A.M. Mathai}\\
Department of Mathematics and Statistics, McGill University\\
Montreal, Canada H3A 2K6, and\\
Centre for Mathematical Sciences Pala Campus\\
Arunapuram P.O., Pala, Kerala 686 574, India\\[0.5cm]
{\bf H.J. Haubold}\\
Office for Outer Space Affairs, United Nations\\
Vienna International Centre, P.O. Box 500, A-1400 Vienna, Austria\\
\end{center}
\noindent
{\bf Abstract.} This paper deals with the solution of unified fractional reaction-diffusion systems. The results are obtained in compact and elegant forms in terms of Mittag-Leffler functions and generalized Mittag-Leffler functions, which are suitable for numerical computation. On account of the most general character of the derived results, numerous results on fractional reaction, fractional diffusion, and fractional reaction-diffusion problems scattered in the literature, including the recently derived results by the authors for reaction-diffusion models, follow as special cases. 

\section{Introduction}

Reaction-diffusion models have found numerous applications in pattern formation in biology, chemistry, and physics.. These systems show that diffusion can produce the spontaneous formation of spatio-temporal patterns.  A piecewise linear approach in connection with the diffusive processes has been   developed by Strier et al. (1995). A similar approach is recently used by Manne et al. (2000) to investigate effects on the propagation of nonlinear wave fronts.\\

The simplest reaction-diffusion model can be described by an equation
\begin{equation}
\frac{\partial N}{\partial t}= D\frac{\partial^2N}{\partial x^2}+\lambda F(N),\;\; N=N(x,t) 
\end{equation}                                                
where D is the diffusion constant and F(N) is a nonlinear function representing reaction kinetics. It is interesting to observe that for $F(N)=\gamma N(1-N),$ (1) reduces to the Fisher-Kolmogorov equation and if, we set $F(N)=\gamma N(1-N^2)$, it reduces to the real Ginsburg-Landau equation. 

A generalization of (1) has been considered by Manne et al. (2000) in the form
\begin{equation}
\frac{\partial^2 N}{\partial t^2}+a\frac{\partial N}{\partial x}=\nu^2\frac{\partial^2 N}{\partial x^2}+\xi^2N(x,t),
\end{equation}
where $\xi$ indicates the strength of the nonlinearity of the system.

In this article, we present a straightforward method for the systematic derivation of the solution of nonlinear reaction-diffusion equations connected with nonlinear waves, which is a generalization of the eq. (2). The results are derived in a closed form by the application of Laplace and Fourier transforms, which are suitable for numerical computation. The present study is a continuation of our investigations reported earlier in the articles of Saxena et al. (2002, 2004, 2004a, 2004b, 2006a, 2006b, 2006c).

\section{Auxiliary Results}

A generalization of the Mittag-Leffler (1903, 1905) function
\begin{equation}                  
E_\alpha(z)=\sum^\infty_{n=0}\frac{z^n}{\Gamma(n\alpha +1)},\;\;(\alpha \in C, Re(\alpha) >0)
\end{equation}                                         
was introduced by Wiman (1905) in the general form 
\begin{equation}                                       
E_{\alpha, \beta}(z)=\sum^\infty_{n=0}\frac{z^n}{\Gamma(n\alpha+\beta)}, (\alpha, \beta \in C, Re (\alpha)>0, Re (\beta) >0).
\end{equation}
The main results of these functions are available in the handbook of Erd\'{e}'lyi, et al. (1955, Section 18.1) and the monographs by Dzherbashyan (1966, 1993). Prabhakar (1971) introduced a generalization of (4) in the form 
\begin{equation}
E^\gamma_{\alpha, \beta}(z)=\sum^\infty_{n=0}\frac{(\gamma)_n z^n}{\Gamma(n\alpha+\beta)(n)!},\;\;(\alpha, \beta,\gamma \in C; Re(\alpha), Re(\beta), Re(\gamma)>0),
\end{equation}                                                              
where $(\gamma)_r$ is the Pochhammer symbol, defined by 
$$(\gamma)_0=1, (\gamma)_r=\gamma(\gamma+1)(\gamma+2)\ldots (\gamma+r-1), r= 1,2,\ldots, \gamma \neq 0.$$                                   
It is an entire function with $\rho=[Re(\alpha)]^{-1}$ (Prabhakar, 1971). Solution of generalized Volterra-type differintegral equations associated with this function as a kernel is derived by Kilbas et al. (2002). A general theory of generalized fractional calculus based on this function has been developed by Kilbas et al. (2004), generalizing the results for Riemann-Liouville fractional integrals and derivatives, which form the backbone of fractional differential equations.\\
For $\gamma=1$, this function coincides with (4), while for $\beta=\gamma=1$ with (3):
\begin{equation}  
E^1_{\alpha, \beta}=E_{\alpha, \beta}(z),\;\;E_{\alpha,1}^1 (z)=E_\alpha (z).
\end{equation}
We also have 
\begin{equation}
\Phi(\alpha, \beta; z)= _1F_1(\alpha; \beta;z)=\Gamma(\beta)E^\alpha_{1,\beta}(z),
\end{equation}   
where $\Phi(\alpha, \beta; z)$ is Kummer's confluent hypergeometric function defined in\\
Erd\'{e}lyi et al. (1953, p. 248, eq. (1)). Prabhakar (1971, p. 8, eq. (2.5)) has shown that
\begin{equation}
L\left\{t^{\gamma-1}e^{-st}E^\delta_{\beta,\gamma}(\omega t^\beta); s\right\}= s^{-\gamma}(1-\omega s^{-\beta})^{-\delta},
\end{equation}
where $Re(\beta)>0, Re(\gamma)>0, Re(s)>0\; \mbox{and}\; s>|\omega|^{\frac{1}{Re(\beta)}}.$

In terms of the H-function, (5) can be expressed in the form 
\begin{equation}
E^\gamma_{\alpha, \beta}(-z)=\frac{1}{\Gamma(\gamma)}H^{1,1}_{1,2}\left[z|^{(1-\gamma, 1)}_{(0,1),(1-\beta,\alpha)}\right],
\end{equation}                                                                   
where $Re(\alpha)>0, Re(\beta)>0, Re(\gamma)>0.$

We will need the definitions of the Laplace and Fourier transforms of a function $N(x,t)$ and their inverses, which are given below:\\
The Laplace transform of a function $N(x,t)$ with respect to $t$ is defined by 
\begin{equation}     
L\left\{N(x,t)\right\}= N^\sim(x,s)=\int^\infty_0 e^{-st}N(x,t)dt, \;\;(t>0),(x\in R),
\end{equation}
where $Re(s) > 0$, and its inverse transform with respect to $s$ is given by
\begin{equation}
L^{-1}\left\{N^\sim(x,s)\right\}=L^{-1}\left\{N^\sim(x,s); t\right\} = \frac{1}{2\pi}\int^{\gamma+ i\infty}_{\gamma-i \infty} e^{st}N^\sim(x,s)ds,
\end{equation}                   
$\gamma$ being a fixed real number. 

The Fourier transform of a function $N(x,t)$ with respect to $x$ is defined by 
\begin{equation}
F\left\{N(x,t)\right\}= F^*(k,t)=\int^\infty_{-\infty} e^{ikx}N(x,t) dx.
\end{equation}       
The inverse Fourier transform with respect to $k$ is given by the formula
\begin{equation}
F^{-1}\left\{F^*(x,t)\right\}=\frac{1}{2\pi}\int^\infty_{-\infty} e^{-ikx}F^*(k,t)dk.
\end{equation}                                                           
The space of functions for which the transforms defined by (10) and (12) exist are denoted by $LF=L(R_+)\times F(R).$                                                                    
In view of the results (Mathai and Saxena, 1978, p. 49), also see (Prudnikov et al., 1989, p. 355, eq. 2.25.3.2), the cosine transform of the H-function is given by 

\begin{eqnarray} 
&&\int^\infty_0 t^{\rho-1}cos(kt)H^{m,n}_{p,q}\left[at^\mu|^{(a_p, A_p)}_{(b_q, B_q)}\right]dt\nonumber\\ 
&=&\frac{\pi}{k^\rho}H^{n+1,m}_{q+1,p+2}\left[\frac{k^\mu}{a}\left|^{(1-b_q, B_q),(\frac{1+\rho}{2}, \frac{\mu}{2})}_{(\rho, \mu),(1-a_p,A_p),(\frac{1+\rho}{2}, \frac{\mu}{2})}\right.\right],
\end{eqnarray}
where\\ 
$Re\left[\rho+\mu Re[\rho+\mu^{min}_{1\leq j\leq m} (\frac{b_j}{B_j})]\right]>1; k, \mu>0; Re\left[\rho+\mu^{max}_{1\leq j\leq n}\left(\frac{a_j-1}{A_j}\right)\right]< \frac{3}{2}.$\\
$|arg\; \alpha|<\frac{1}{2}\pi \theta; \theta=\sum^n_{j=1}A_j-\sum^p_{j=n+1}A_j+\sum^m_{j=1}B_j-\sum^q_{j=m+1}B_j>0.$
     
The Riemann-Liouville fractional integral of order $\nu$ is defined by (Miller and Ross, 1993, p. 45) 
\begin{equation}    
_0D_t^{-\nu}N(x,t)=\frac{1}{\Gamma(\nu)}\int_0^t(t-u)^{\nu-1}N(x,u)du,
\end{equation}
where $Re(\nu)>0.$

Following Samko et al. (1990, p. 37), we define the fractional derivative of order $\alpha>0$ in the form
\begin{equation}
_0D_t^\alpha N(x,t)=\frac{1}{\Gamma(n-\alpha)}\frac{d^n}{dt^n}\int^t_0\frac{N(x,u)du}{(t-u)^{\alpha-n+1}},\;\;(t>0), (n=[\alpha]+1),
\end{equation}
where $[\alpha]$ means the integral part of the number $\alpha$.

From Erd\'{e}lyi et al. (1954b, p. 182), we have 
\begin{equation}
L\left\{_0D_t^{-\alpha}N(x,t)\right\}=s^{-\alpha}N^\sim(x,s),
\end{equation}
where $Re(s)> 0$ and $Re(\alpha)>0$.                                                                    
     
The Laplace transform of the fractional derivative, defined by (16), is given by Oldham and Spanier, (1974, p. 134, eq. (8.1.3))
\begin{equation}
L\left\{_0D_t^\alpha N(x,t)\right\}=s^\alpha N^\sim(x,s)-\sum^n_{r=1}s^{r-1}\; _0D_t^{\alpha-r}N(x,t)|_{t=0},\;\;n-1<\alpha\leq n.
\end{equation}
In certain boundary-value problems arising in the theory of viscoelasticity and in the hereditary solid mechanics, the following fractional derivative of order $\alpha>0$ is introduced by Caputo (1969) in the form
\begin{eqnarray}
_0D_t^\alpha f(x,t)& = & \frac{1}{\Gamma(m-\alpha)}\int_0^t\frac{f^{(m)}(x,\tau)d\tau}{(t-\tau)^{\alpha+1-m}}\nonumber\\
&&m-1<\alpha\leq m, \;Re(\alpha)>0,\;m\in \mbox{N}\nonumber\\
&=&\frac{\partial^mf(x,t)}{\partial t^m},\; \mbox{if}\; \alpha = m.
\end{eqnarray}  
Caputo (1969) has given the Laplace transform of the derivative as 
\begin{equation}
L\left\{_0D_t^\alpha f(x,t)\right\}= s^\alpha f^\sim (x,s) - \sum ^{m-1}_{r=0}s^{\alpha-r-1}f^{(r)}(x,0), \; (m-1<\alpha\leq m),
\end{equation}
where $f^\sim(x,s)$ is the Laplace transform of  $f(x,t)$.                                                 
     
The above formula is useful in deriving the solution of differintegral equations of fractional order governing certain physical problems of reaction and diffusion. 
 
We also need the Weyl fractional operator, defined by
\begin{equation} 
_{-\infty}D^\mu_xf(x,t)= \frac{1}{\Gamma(n-\mu)}\frac{d^n}{dt^n}\int^t_{-\infty}\frac{f(x,u)du}{(t-u)^{\mu-n+1}},
\end{equation}    
where $n=[\mu]+1$ is an integral part of $\mu>0$.
     
The Fourier transform of the Caputo derivative is given by (Metzler and Klafter, 2000, p. 59, A.11)
\begin{equation}
F\left\{_{-\infty}D^\mu_xf(x,t)\right\}= (ik)^\mu f^*(k,t).
\end{equation}                                                                                                     
Following the convention initiated by Compte (1996), we suppress the imaginary unit in Fourier space by adopting a slightly modified form of the result (21) in our investigations (Metzler and Klafter, 2000, p. 59, A.12):
\begin{equation}
F\left\{_{-\infty}D_x^\mu f(x,t)\right\}=-|k|^\mu f^*(k,t).
\end{equation}                                                              

\section{Inverse Laplace Transform of Certain\\
Algebraic Functions}

This section deals with the evaluation of the inverse Laplace transforms of certain algebraic functions which are directly applicable in the analysis that follows.

It will be shown here that
\begin{eqnarray}
&&L^{-1}\left\{\frac{s^{\rho-1}}{s^\alpha + as^\beta+bs^\gamma+c};t\right\}=\sum^\infty_{r=0}(-1)^r\nonumber\\
&&\sum^r_{l=0}(^r_l)a^l b^{r-l}t^{(\alpha-\gamma)r+(\gamma-\beta)l+\alpha-\rho}E^{r+1}_{\alpha,(\alpha-\gamma)r+(\gamma-\beta)l+\alpha-\rho+1}(-ct^\alpha),
\end{eqnarray}                                                                                                     
where\\
$Re(\alpha)>0, Re(\beta)>0, Re(\gamma)>0, Re(s)>0, Re(\rho)>0, \left|\frac{as^\beta+bs^\gamma}{s^\alpha +c}\right|<1,\;\;E^\delta_{\beta, \gamma}(.)$  is the generalized Mittag-Leffler function defined by eq. (5) and provided that the series in eq. (24) is convergent.\\
{\bf Proof.}  Assume that $\alpha>\gamma>\beta$. We have

\begin{eqnarray}
&&\frac{s^{\rho-1}}{s^\alpha+as^\beta+bs^\gamma+c}\nonumber\\
&=&\frac{s^{\rho-\beta-1}}{(s^{\alpha-\beta}+cs^{-\beta})(1+\frac{a+bs^{\gamma-\beta}}{s^{\alpha-\beta}+cs^{-\beta}})}\nonumber\\
&=&s^{\rho-\beta-1}\sum^\infty_{r=0}\frac{(-1)^r(a+bs^{\gamma -\beta})^r}{(s^{\alpha-\beta}+cs^{-\beta})^{r+1}}\nonumber\\
&=&\sum^\infty_{r=0}(-1)^r\sum^r_{l=0}(^r_l)a^lb^{r-l}\frac{s^{\rho+(\gamma-\alpha)r+(\beta-\gamma)l-\alpha-1}}{(1+cs^{-\alpha})^{r+1}}.
\end{eqnarray} 
On term by term inverting (25) with the help of the formula (8), it readily gives the desired result. The term by term inversion is justified by virtue of a theorem due to Doetsch (1956, \S 22).
 
(i) From (25), it readily follows that      
\begin{equation}
L^{-1}\left\{\frac{s^{\rho-1}}{s^\alpha+as^\beta+b};t\right\}=t^{\alpha-\rho}\sum^\infty_{r=0}(-a)^r t^{(\alpha-\beta)r}E^{r+1}_{\alpha,\alpha+(\alpha-\beta)r-\rho+1}(-bt^\alpha),
\end{equation}
where $Re(\alpha)>0, Re(\beta)>0, Re(\rho)>0, Re(s)>0, \left|\frac{as^\beta}{s^\alpha+b}\right|<1$  and provided that the series in (26) is convergent.
 
Some special cases of the result (26) are worth mentioning.

For $\rho=\alpha$, (26) reduces to the following result given by the authors (Saxena et al., 2006b):
\begin{equation}
L^{-1}\left\{\frac{s^{\alpha-1}}{s^\alpha+as^\beta+b}; t\right\}= \sum^\infty_{r=0}(-a)^r t^{(\alpha-\beta)r}E^{r+1}_{\alpha,(\alpha-\beta)r+1} (-bt^\alpha), 
\end{equation}
where $Re(\alpha)>0, Re(\beta)>0, Re(s) >0, |\frac{as^\beta}{s^\alpha +b}|<1$ and provided that the series in (27) is convergent.
 
Further, if we set $\rho=\beta$, then we arrive at another result given by the authors: 
\begin{equation}
L^{-1}\left\{\frac{s^{\beta-1}}{s^\alpha+as^\beta+b}; t\right\}= t^{\alpha-\beta}\sum^\infty_{r=0}(-a)^r t^{(\alpha-\beta)r}E^{r+1}_{\alpha,(\alpha-\beta)(r+1)+1}(-bt^\alpha),
\end{equation}          
where $Re(\alpha)>0, Re(\beta)>0, Re(s)>0,\left|\frac{as^\beta}{s^\alpha+b}\right|<1, \alpha>\beta$, and provided that the series in (28) is convergent.
 
(ii) If we set $\rho=1$ in (24), we obtain the following result:
\begin{eqnarray}
&&L^{-1}\left\{\frac{1}{s^\alpha+as^\beta+bs^\gamma+c};t\right\}= t^{\alpha-1}\sum^\infty_{r=0}(-1)^r\nonumber\\
&\times&  \sum^r_{l=0}(^r_l)a^l b^{r-l}t^{(\alpha-\gamma)r+(\gamma-\beta)l+\alpha-1}\nonumber\\
&\times&E^{r+1}_{\alpha, \alpha+(\alpha-\beta)r+(\gamma-\beta)l+\alpha}(-ct^\alpha),
\end{eqnarray} 
where $Re(\alpha)>0, Re(\beta)>0, Re(\gamma)>0, Re(s)>0, \left|\frac{as^\beta+bs^\gamma}{s^\alpha+c}\right|<1,$ and provided that the series in (29) is convergent.

(iii) Next, if we set $\rho=\alpha$ in (24), we then arrive at the following result
\begin{eqnarray}
&& L^{-1}\left\{\frac{s^{\alpha-1}}{s^\alpha+as^\beta+bs^\gamma+c}; t\right\}\nonumber\\
&=& \sum^\infty_{r=o}(-1)^r\sum^r_{l=0}(^r_l)a^l b^{r-l}t^{(\alpha-\gamma)r+(\gamma-\beta)l}E^{r+1}_{\alpha,(\alpha-\gamma)r+(\gamma-\beta)l+1}(-ct^\alpha),
\end{eqnarray}
where $Re(\alpha)>0, Re(\beta)>0, Re(\gamma)>0, Re(s)>0, \left|\frac{as^\beta+bs^\gamma}{s^\alpha+c}\right|<1$, and provided that the series in (30) is convergent.

(iv) Similarly, $\rho=\beta$  in (24) will yield the formula
\begin{eqnarray}
 && L^{-1}\left\{\frac{t^{\beta-1}}{s^\alpha+as^\beta+bs^\gamma+c}; t\right\}\nonumber\\
&=& \sum^\infty_{r=0}(-1)^r\sum^r_{l=0}(^r_l)a^l b^{r-l}t^{(\alpha-\gamma)r+(\gamma-\beta)l+\alpha-\beta}E^{r+1}_{\alpha,(\alpha-\gamma)r+(\gamma-\beta)l+\alpha-\beta+1}(-ct^\alpha)\nonumber\\
&&             
\end{eqnarray}
where $Re(\alpha)>0, Re(\beta)>0, Re(\gamma)>0, Re(s)>0, \left|\frac{as^\beta+bs^\gamma}{s^\alpha+c}\right|<1$, and provided that the series in (31) is convergent.

(v) When $\rho=\gamma$,(24) readily gives the formula
\begin{eqnarray}
&& L^{-1}\left\{\frac{t^{\gamma-1}}{s^\alpha+as^\beta+bs^\gamma+c}; t\right\}\nonumber\\
&=& \sum^\infty_{r=o}(-1)^r\sum^r_{l=0}(^r_l)a^l b^{r-l}t^{(\alpha-\gamma)r+(\gamma-\beta)l+\alpha-\gamma}E^{r+1}_{\alpha,(\alpha-\gamma)r+(\gamma-\beta)l+\alpha-\gamma+1}(-ct^\alpha),\nonumber\\ 
&&             
\end{eqnarray}
where $Re(\alpha)>0, Re(\beta)>0, Re(\gamma)>0, Re(s)>0, \left|\frac{as^\beta+bs^\gamma}{s^\alpha +c}\right|<1$, and provided that the series in (32) is convergent. 

(vi) When $\rho=\alpha+\beta$ then (24) reduces to
\begin{eqnarray}
&&L^{-1}\left\{\frac{s^{\alpha+\beta-1}}{s^\alpha+as^\beta+bp^\gamma+c}; t\right\}\nonumber\\
&=&\sum^\infty_{r=0}(-1)^r\sum^r_{l=0}(^r_l)a^l b^{r-l}t^{(\alpha-\gamma)r+(\gamma-\beta)l-\beta}E^{r+1}_{\alpha,(\alpha-\gamma)r+(\gamma-\beta)l-\beta+1}(-ct^\alpha),\nonumber\\
&&
\end{eqnarray}  
where $Re(\alpha)>0, Re(\beta)>0, Re(\gamma)>0, Re(\alpha+\beta)>0, \left|\frac{as^\beta+bs^\gamma}{s^\alpha +c}\right|<1$, and provided that the series in (33) is convergent. 

\section{A General Case}

In this Section, it will be shown that if $Re(\alpha_j)>0\;\; (j=1,\ldots,n), Re(s)>0$, then
\begin{eqnarray}
&&L^{-1}\left[\frac{s^{\rho-1}}{a_0+a_1 s^{\alpha_1}+a_2s^{\alpha_2}+a_{n-2}s^{\alpha_3}+\ldots+a_ns^{\alpha_n}+a_{n+1}s^{\alpha_{n+1}}};t\right]\nonumber\\
&=&\frac{1}{a_1}\sum^\infty_{m=0}(-1)^m\sum_{^{r_1+r_2+\ldots+r_n=m}_{ r_1\geq 0\ldots r_n\geq 0}}\frac{m!}{r_1!\ldots r_n!}\left\{\Pi^n_{j=1}\left(\frac{a_{j+1}}{a_1}\right)^{r_j}\right\}t^{A-\rho}E^{m+1}_{\alpha_1, A-\rho+1}\nonumber\\
&\times&\left(-\frac{a_0}{a_1}t^{\alpha_1}\right), a_1\neq 0
\end {eqnarray}  
where $A = \sum^n_{j=1}(\alpha_2-\alpha_{j+1})r_j+\alpha_1(1+m)-\alpha_2m, \left|\frac{\sum^n_{r=1}\left\{(a_{r+1}/a_1)s^{\alpha_{r+1}-\alpha_2}\right\}}{(s^{\alpha_1-\alpha_2}+(a_0/a_1)s^{-\alpha_2})}\right|<1,$
provided that the series in (34) is convergent.\\      
{\bf Proof:}  Let us assume that $\alpha_2>\alpha_{j+1}, (j=2,3,\ldots,n)$. We have 
\begin{eqnarray}
&&\frac{s^{\rho-1}}{a_0+a_1 s^{\alpha_1}+a_2 s^{\alpha_2}+a_3 s^{\alpha_3}+\ldots+a_n s^{\alpha_n}+a_{n+1}s^{n+1}} = \nonumber\\
&=&\frac{a_1^{-1}s^{\rho-\alpha_2-1}}{[s^{\alpha_1-\alpha_2}+(a_0/a_1)s^{-\alpha_2}]\left[1+\frac{\sum^n_{r=1}\left\{(a_{r+1}/a_1)s^{\alpha_{r+1}-\alpha_2}\right\}}{(s^{\alpha_1-\alpha_2}+(a_0/a_1)s^{-\alpha_2})}\right]}\nonumber\\
&=& a_1^{-1}s^{\rho-\alpha_2-1}\sum^\infty_{m=0}(-1)^m\frac{\left[\sum^n_{r=1}\left\{(a_{r+1}/a_1)s^{\alpha_{r+1}-\alpha_2}\right\}\right]^m}{(s^{\alpha_1-\alpha_2}+(a_0/a_1)s^{-\alpha_2})^{m+1}}.
\end {eqnarray}  
On applying the multinomial theorem (Abramowitz and Stegun, 1968, p. 823, para 24.1.2):
$$(x_1+x_2+\ldots +x_m)^n=\sum(n;n_1,n_2,\ldots,n_m)x_1^{n_1}x_2^{n_2}\ldots x_m^{n_m},$$
summed over $n_1+n_2+\ldots+n_m = n$; the above line transforms into the form
\begin{eqnarray}
&=&\frac{1}{a_1}\sum^\infty_{m=0}(-1)^m  \sum_{^{r_1+\ldots+r_n=m}_{r_1\geq 0,\ldots,r_n\geq 0}}\frac{m!}{r_1!\ldots r_n!}\left\{\Pi^n_{j=1}\left(\frac{a_{j+1}}{a_1}\right)^{r_j}\right\}\nonumber\\
&\times&\frac{s^{\sum^n_{j=1}(\alpha_{j+1}-\alpha_2)r_j+\rho+\alpha_2m-(1+m)\alpha_1-1}}{(1+(a_0/a_1)s^{-\alpha_1})^{m+1}}.
\end{eqnarray}
Interpreting the above expression with the help of the formula (8), we obtain the desired result. This completes the Proof of (34). Some interesting special cases of the general result (34) are given below. More results can be derived by specializing the parameters and the variable and applying certain known theorems of Laplace and inverse transforms.

(i) If we set $\rho=\alpha_1$, then (34) reduces to
\begin{eqnarray}
&&L^{-1}\left[\frac{s^{\alpha_1-1}}{a_0+a_1s^{\alpha_1}+a_2s^{\alpha_2}+a_{n-2}s^{\alpha_3}+\ldots+a_n s^{\alpha_n}+a_{n+1}s^{\alpha_{n+1}}};t\right]\nonumber\\
&=& a_1^{-1}\sum^\infty_{m=0}(-1)^m \sum_{^{r_1+r_2+\ldots +r_n=m}_{r\geq 0\ldots r_n\geq 0}}\frac{m!}{r_1!\ldots r_n!}\left\{\Pi^n_{j=1}\left(\frac{a_{j+1}}{a_1}\right)^{r_j}\right\}t^{A-\alpha_1}\nonumber\\
&\times& E^{m+1}_{\alpha_1,A-\alpha_1+1}\left(-\frac{a_0}{a_1}t^{\alpha_1}\right), a_1\neq 0,
\end{eqnarray}           
where\\ 
$Re(\alpha_j)>0\;\; (j=1,\ldots, n), Re(s)>0, \left|\frac{\sum^n_{r=1}\left\{(a_{r+1}/a_1)s^{\alpha_{r+1}-\alpha_2}\right\}}{(s^{\alpha_1-\alpha_2}+(a_0/a_1)s^{-a_2})}\right|<1, A=\sum^n_{j=1}(\alpha_2-\alpha_{j+1})r_j+\alpha_1(1+m)-\alpha_2m,$
and provided that the series in (37) is convergent.

(ii) On the other hand, if we take $\rho=\alpha_2$ in (34), it gives the formula
      
\begin{eqnarray}
&&L^{-1}\left[\frac{s^{\alpha_{2-1}}}{a_0+a_1s^{\alpha_1}+a_2s^{\alpha_2}+a_{n-2}s^{\alpha_3}+\ldots+a_n s^{\alpha_n}+a_{n+1}s^{\alpha_{n+1}}};t\right]\nonumber\\
&=& a_1^{-1}\sum^\infty_{m=0}(-1)^m \sum_{^{r_1+r_2+\ldots +r_n=m}_{r_1\geq 0\ldots r_n\geq 0}}\frac{m!}{r_1!\ldots r_n!}\left\{\Pi^n_{j=1}\left(\frac{a_{j+1}}{a_1}\right)^{r_j}\right\}t^{A-\alpha_2}\nonumber\\           
&\times& E^{m+1}_{\alpha_1,A-\alpha_2+1}\left(-\frac{a_0}{a_1}t^{\alpha_1}\right), a_1\neq 0,
\end{eqnarray}                                                             
where $Re(\alpha_j)>0\;\; (j=1,\ldots, n), Re(s)>0,\left(\left|\frac{\sum^n_{r=1}\left\{(a_{r+1}/a_1)s^{\alpha_{r+1}-\alpha_2}\right\}}{(s^{\alpha_1-\alpha_2}+(a_0/a_1)s^{-a_2})}\right|\right)<1,$ \\
$A=\sum^n_{j=1}(\alpha_2-\alpha_{j+1})r_j+\alpha_1(1+m)-\alpha_2 m,$,
and provided that the series in (38) is convergent .

(iii) Next, if we set $\rho=1$ in (34), then we arrive at the following inversion formula given earlier by Podlubny (1999) in a slightly different form:
\begin{eqnarray}
&&L^{-1}\left[\frac{1}{a_0+a_1s^{\alpha_1}+a_2s^{\alpha_2}+a_{n-2}s^{\alpha_3}+\ldots+a_n s^{\alpha_n}+a_{n+1}s^{\alpha_{n+1}}};t\right]\nonumber\\
&=& a_1^{-1}\sum^\infty_{m=0}(-1)^m \sum_{^{r_1+r_2+\ldots +r_n=m}_{r\geq 0\ldots r_n\geq 0}}\frac{m!}{r_1!\ldots r_n!}\left\{\Pi^n_{j=1}\left(\frac{a_{j+1}}{a_1}\right)^{r_j}\right\}t^{A-1}\nonumber\\           
&\times& E^{m+1}_{\alpha_1,A}\left(-\frac{a_0}{a_1}t^{\alpha_1}\right), a_1\neq 0,
\end{eqnarray}           
where
$$Re(\alpha_j)>0(j=1,\ldots,n), Re(s)>0, \left|\frac{\sum^n_{r=1}\left\{(a_{r+1}/a_1)s^{\alpha_{r+1}-\alpha_2}\right\}}{(s^{\alpha_1-\alpha_2}+(a_0/a_1)s^{-\alpha_2})}\right|<1,$$
$A= \sum^n_{j=1}(\alpha_2-\alpha_{j+1})r_j+\alpha_1(1+m)-\alpha_2m$, and provided that the series in (4.5) is convergent.      
 
(iv) When $\rho=\alpha_1+\alpha_2,$ (34), then it reduces to the formula 
\begin{eqnarray}
&&L^{-1}\left[\frac{s^{\alpha_1+\alpha_2-1}}{a_0+a_1s^{\alpha_1}+a_2s^{\alpha_2}+a_{n-2}s^{\alpha_3}+\ldots+a_n s^{\alpha_n}+a_{n+1}s^{\alpha_{n+1}}};t\right]\nonumber\\
&=& a_1^{-1}\sum^\infty_{m=0}(-1)^m \sum_{^{r_1+r_2+\ldots +r_n=m}_{r\geq 0\ldots r_n\geq 0}}\frac{m!}{r_1!\ldots r_n!}\left\{\Pi^n_{j=1}\left(\frac{a_{j+1}}{a_1}\right)^{r_j}\right\}t^{A_1}\nonumber\\           
&\times& E^{m+1}_{\alpha_1,A_1+1}\left(-\frac{a_0}{a_1}t^{\alpha_1}\right), a_1\neq 0,
\end{eqnarray}   
where $Re(\alpha_j)>0. (j=1,2,\ldots, n): \alpha_2>\alpha_{j+1}(j=2,3,\ldots, n), \alpha_1 \neq 0, Re(s)>0,$
$$\left|\frac{\sum^n_{r=1}\left\{(a_{r+1}/a_1)s^{\alpha_{r+1}-\alpha_2}\right\}}{(s^{\alpha_1-\alpha_2}+(a_0/a_1)s^{-\alpha_2})}\right|<1, \;\;A_1=\sum^n_{j=1}(\alpha_2-\alpha_{j+1}r_j+\alpha_1 m-\alpha_2(m+1),$$
and provided that the series in (40) is convergent.  
In what follows, $e^\delta_{\beta, \gamma}(.)$ will be employed to denote the generalized Mittag-Leffler function, defined by (5).

\section{Solution of Fractional Reaction-Diffusion\\
Equation}
        
In this Section, it is proposed to derive the solution of the fractional diffusion system connected with nonlinear waves governed by the eq. (41). This system is a generalized form of the reaction-diffusion equation recently studied by Manne et al. (2000). The result is given in the form of the following \\
{\bf Theorem 1.} Consider the fractional reaction-diffusion equation associated with Caputo derivatives in the form 
\begin{eqnarray}
&&_0D_t^\alpha N(x,t)+a_0\; D^\beta_t N(x,t)+b_0D_t^\gamma N(x,t)\nonumber\\
&=&\nu^2\;_{-\infty}D_x^\eta N(x,t)+\xi^2N(x,t)+\phi(x,t)\\
&&(0\leq\alpha\leq 1,0\leq \beta \leq 1,0 \leq \gamma \leq 1, \eta > 0)\nonumber
\end{eqnarray}                                                                             
with initial conditions 
\begin{equation}
N(x,0)=f(x), \;x\in  \Re, lim_{x\rightarrow \pm\infty} N(x,t)=0, t>0,
\end{equation}
where $\nu^2$  is a diffusion constant, $\xi$ is a constant which describes the nonlinearity in the system, and $\varphi(x,t)$ is nonlinear function which belongs to the area of reaction-diffusion, then there holds the following formula for the solution of (41):
\begin{eqnarray}
&&N(x,t)=\frac{1}{2\pi}\sum^\infty_{r=0}(-1)^r \int^\infty_{-\infty} e^{-ikx}\sum^r_{l=0}(^r_l)a^l b^{r-l}t^{(\alpha-\gamma)r+(\gamma-\beta)l}f^*(k)\nonumber\\
&&\left\{E^{r+1}_{\alpha,(\alpha-\gamma)r+(\gamma-\beta)l+1}(-ct^\alpha)\right.\nonumber\\
&+&\left. at^{\alpha-\beta} E^{r+1}_{\alpha,(\alpha-\gamma)r+(\gamma-\beta)l+1}(-ct^\alpha)+b t^{\alpha-\gamma}
E^{r+1}_{\alpha,(\alpha-\gamma)r+(\gamma-\beta)l+1}(-ct^\alpha)\right\}dk\nonumber\\
&+&\sum^\infty_{r=0}\frac{(-1)^r}{2\pi}\sum^r_{l=0}(^r_l)a^l b^{r-l}\int^t_0 \xi^{(\alpha-\gamma)r+(\gamma-\beta)l+\alpha-1}d\xi \int^\infty_{-\infty}e^{-ikx}\varphi^*(k,t-\xi)\nonumber\\
&&E^{r+1}_{\alpha,(\alpha-\gamma)r+(\gamma-\beta)l+\alpha} (-b\xi^\alpha)dkd\xi
\end{eqnarray}
where $c=v^2|k|^\delta-\xi^2$ and provided that the series and integrals in (43) are convergent.\\
{\bf Proof.}  Applying the Laplace transform with respect to the time variable
 $t$ and using the boundary conditions, we find that 
\begin{eqnarray}
&&s^\alpha N^\sim(x,s)- s^{\alpha-1}f(x) + as^\beta N^\sim(x,s)-as^{\beta-1} f(x)+\nonumber\\
&& bs^\gamma N^\sim (x,s)-bs^{\gamma-1}f(x)\nonumber\\             
&=&\nu^2\;_{-\infty}D_x^\eta N^\sim(x,s)+\xi^2N^\sim(x,s)+\varphi^\sim (x,s)
\end{eqnarray}
If we apply the Fourier transform with respect to the space variable $x$, it yields           
\begin{eqnarray}
&&s^\alpha N^{\sim *}(k,s)- s^{\alpha-1}f^*(k) + as^\beta N^{\sim *}(k,s), -as^{\beta-1} f^*(k)+\nonumber\\
&&bs^\gamma N^{\sim *} (k,s)-bs^{\gamma-1}f^*(k)\nonumber\\             
&=&-\nu^2|k|^\eta N^{\tilde *}(k,s)+ \xi^2 N^{\sim *}(k,s)+\varphi^{\sim *} (k,s).
\end{eqnarray}
Solving for $N^{\sim *}(k,s)$, it gives
\begin{equation}
N^{\sim *} (k,s)=\frac{(s^{\alpha-1}+as^{\beta-1}+bs^{\gamma-1})f^*(k)+\varphi^{\sim *} (k,s)}{s^\alpha+as^\beta+bs^\gamma+c}
\end{equation}                            
where  $c=\nu^2|k|^\eta - \xi^2.$

To invert the eq. (46), it is convenient to first invert the Laplace transform and then the Fourier transform. Inverting the Laplace transform with the help of the results (29), (30), (31) and (32), it yields
\begin{eqnarray}
&&N^*(k,t)=\sum^\infty_{r=0}(-1)^r \sum^r_{l=0}(^r_l)a^l b^{r-l} t^{(\alpha-\gamma)r+(\gamma-\beta)l}f^*(k)\nonumber\\
&&\left[E^{r+1}_{\alpha,(\alpha-\gamma)r+(\gamma-\beta)l+1}(-ct^\alpha)\right.\nonumber\\
&+&\left. at^{\alpha-\beta}E^{r+1}_{\alpha,(\alpha-\gamma)r+(\gamma-\beta)l+\alpha-\beta+1}(-ct^\alpha)+\right.\nonumber\\
&&\left. bt^{\alpha-\gamma}E^{r+1}_{\alpha,(\alpha-\gamma)r+(\gamma-\beta)l+\alpha-\gamma+1}(-ct^\alpha)\right]\nonumber\\    
&+&\sum^\infty_{r=0}(-1)^r\sum^r_{l=0}(^r_l)a^l b^{r-l}\nonumber\\
&&\int_0^t\varphi^*(k,t-\xi)\xi^{(\alpha-\gamma)r+(\gamma-\beta)l+\alpha-1}E^{r+1}_{\alpha,(\alpha-\gamma)r+(\gamma-\beta)l+\alpha}(-c\xi^\alpha)d\xi.\nonumber\\
&&
\end{eqnarray}
Finally, the inverse Fourier transform of (47) gives the desired solution in the form 
\begin{eqnarray}
&&N(x,t)=\frac{1}{2\pi}\sum^\infty_{r=0}(-1)^r\int^\infty_{-\infty} e^{-ikx}\sum^r_{l=0}(^r_l)a^l b^{r-l}t^{(\alpha-\gamma)r+(\gamma-\beta)l}f^*(k)\nonumber\\
&&\left[E^{r+1}_{\alpha,(\alpha-\gamma)r+(\gamma-\beta)k+1}(-ct^\alpha)\right.\nonumber\\
&+& \left.abt^{\alpha-\beta}E^{r+1}_{\alpha,(\alpha-\gamma)r+(\gamma-\beta)k+11}(-ct^\alpha)+ \right.\nonumber\\
&&\left. b t^{\alpha-\gamma} E^{r+1}_{\alpha,(\alpha-\gamma)r+(\gamma-\beta)l+1}(-ct^\alpha)\right]dk\nonumber\\
&+&\sum^\infty_{r=0}\frac{(-1)^r}{2\pi}\sum^r_{l=0}(^r_l)a^l b^{r-l}\int^t_0\xi^{(\alpha-\gamma)r+(\gamma-\beta)l+\alpha-1}d\xi\times \nonumber\\
&&\;\;\times\int^\infty_{-\infty}e^{-ikx}\varphi^*(k,t-\xi)E^{r+1}_{\alpha,(\alpha-\gamma)r+(\gamma-\beta)l+\alpha}(-c\xi^\alpha)dk,
\end{eqnarray}
This completes the proof of the theorem.

\section{Interesting Special Cases and Fundamental\\
Solutions}

Some special cases of the theorem 1 are given below.
(i) If we set $\alpha=\beta=\gamma=\eta=1/2$, the theorem reduces to \\
{\bf Corollary 1.1.} Consider the fractional reaction-diffusion equation 
\begin{eqnarray}
&&_0D_t^{1/2}N(x,t)+a_0\;D_t^{1/2}N(x,t)+b_0D_t^{1/2}N(x,t)\nonumber\\
&=&\nu^2_{\-\infty}D_x^{1/2}N(x,t)+\xi^2N(x,t)+\varphi(x,t),
\end{eqnarray}                                                  
with initial conditions 
\begin{equation}
N(x,0)=f(x), x\in\Re,\;\;lim_{x\rightarrow \pm\infty} N(x,t)=0,t>0,
\end{equation}
where $\nu^2$  is a diffusion constant, $\xi$  is a constant which  describes the nonlinearity in the system, and is nonlinear function which belongs to the area of reaction-diffusion, then there holds the following formula for the solution of (49):
\begin{eqnarray}
N(x,t)&=&\frac{(1+a+b)}{2\pi}\sum^\infty_{r=0}(-1)^r\int^\infty_{-\infty} e^{-ikx}\sum^r_{l=0}(^r_l)a^lb^{r-l}f^*(k)E^{r+1}_{1/2,1}(-c^*t^\alpha)dk\nonumber\\
&+&\sum^\infty_{r=0}\frac{(-1)^r}{2\pi}\sum^r_{l=0}(^r_l)a^lb^{r-l}\int^t_0\int^\infty_{-\infty}e^{-ikx}\varphi^*(k,t-\xi)\xi^{-1/2}\times\nonumber\\
&\times& E^{r+1}_{1/2,1/2}(-c^*\xi^\alpha)dkd\xi
\end{eqnarray}
where $c^*=v^2|k|^{1/2}-\xi^2.$

(ii) When $f(x)=\delta(x),$ where $\delta(x)$ is the Dirac-delta function, the theorem yields the following result \\
{\bf Corollary 1.2.} Consider the fractional reaction-diffusion system governed by the equation
\begin{eqnarray}                
&&_0D_t^\alpha N(x,t)+a_0 D_t^\beta N(x,t)+b_0 D_t^\gamma N(x,t)\nonumber\\
&=&\nu^2 _{-\infty}D^\eta_x N(x,t)+\xi^2 N(x,t)+\varphi(x,t)\\
&&(0\leq \alpha \leq 1,0\leq \beta \leq 1,0 \leq \gamma \leq 1,\eta >0)\nonumber,
\end{eqnarray}
subject to the initial conditions 
\begin{equation} 
N(x,0) = \delta(x), x\in\Re, lim_{x\rightarrow \pm\infty} N(x,t) = 0, t> 0,
\end{equation}
where $\delta(x)$  is the Dirac-delta function. Here $\xi$  is a constant that describes the nonlinearity in the system, and $\varphi(x,t)$ is nonlinear function which belongs to the area of reaction-diffusion. Then there exists the following formula for the fundamental solution of (52) subject to the initial conditions (53):
\begin{eqnarray}
N(x,t)&=&\frac{1}{2\pi}\sum^\infty_{r=0} (-1)^r\int^\infty_{-\infty} e^{-ikx}\sum^r_{l=0}(^r_l)a^l b^{r-l}t^{(\alpha-\gamma)r+(\gamma-\beta)l}\nonumber\\
&&\left\{E^{r+1}_{\alpha,(\alpha-\gamma)r+(\gamma-\beta)l+1}(-ct^\alpha)\right.\nonumber\\
&+& \left. at^{\alpha-\beta}E^{r+1}_{\alpha,(\alpha-\gamma)r+(\gamma-\beta)k+11}(-ct^\alpha)+bt^{\alpha-\gamma}\right.\nonumber\\
&&\left. E^{r+1}_{\alpha,(\alpha-\gamma)r+(\gamma-\beta)l+1}(-ct^\alpha)\right\}dk\nonumber\\
&+&\sum^\infty_{r=0}\frac{(-1)^r}{2\pi}\sum^r_{l=0}(^r_l)a^l b^{r-l}\int^t_0\nonumber\\
&& \int^\infty_{-\infty} e^{-ikx}\varphi^*(k,t-\xi)\xi^{(\alpha-\gamma)r+(\gamma-\beta)l+\alpha-1}\nonumber\\
&&E^{r+1}_{\alpha,(\alpha-\gamma)r+(\gamma-\beta)l+\alpha}(-c\xi^\alpha)dkd\xi
\end{eqnarray}
where $b=\nu^2|k|^\eta - \xi^2$, and provided that the series and integrals in (54) are convergent.
     
As $b\rightarrow 0$, we obtain the following result recently obtained  by the authors
(Saxena et al., 2006b):\\
{\bf Corollary 1.3.} Consider the fractional reaction-diffusion equation 
\begin{eqnarray}                              
&& _0D_t^\alpha N(x,t) + a\; _0 D_t^\beta N(x,t)\nonumber\\
&=& \nu^2 _{-\infty}D_x^\gamma N(x,t) + \xi^2 N(x,t) + \varphi (x,t)\\
&&(0\leq \alpha \leq 1,0 \leq \beta \leq 1, \gamma>0)\nonumber
\end{eqnarray}                       
with initial conditions 
\begin{equation}
N(x,0)=f(x), \in \Re, lim_{x\rightarrow \pm\infty} N(x,t) = 0, t>0,
\end{equation}
where $\nu^2$  is a diffusion constant, $\xi$  is a constant which  describes the nonlinearity in the system, and $\phi(x,t)$ is nonlinear function which belongs to the area of reaction-diffusion, then there holds the following formula for the solution of (55).
\begin{eqnarray}
&& N(x,t)\nonumber\\
&=&\sum^\infty_{r=0}\frac{(-1)^r}{2\pi}\int^\infty_{-\infty} t^{(\alpha-\beta)r}f^*(k)exp(-kx)\left[E^{r+1}_{\alpha,(\alpha-\beta)r+1}(-bt^\alpha)+t^{\alpha-\beta}\right.\nonumber\\
&&\left.E^{r+1}_{\alpha,(\alpha-\beta)(r+1)+1}(-bt^\alpha)\right]dk\nonumber\\
&+&\sum^\infty_{r=0}\frac{(-a)^r}{2\pi}\int^t_0 \xi^{\alpha+(\alpha-\beta)r-1}\int^\infty_{-\infty}\varphi(k,t-\xi)exp(-ikx)E^{r+1}_{\alpha,\alpha+(\alpha-\beta)r}\nonumber\\
&&(-b\xi^\alpha)dkd \xi,
\end{eqnarray}         
where $\alpha>\beta$   and $b=\nu^2|k|^\gamma-\xi^2.$
Now if we set $f(x)=\delta(x), \eta=2$, $\alpha$ is replaced by $2\alpha$ and $\beta$  by $\alpha, \gamma=0$, and use the result (Saxena et al., 2006b, eqs. (24) and (28)),the following result is obtained.\\
{\bf Corollary 1.4} Consider the following reaction-diffusion system
\begin{equation}
\frac{\partial^{2\alpha}N(x,t)}{\partial t^{2\alpha}}+a\frac{\partial^\alpha N(x,t)}{\partial t^\alpha}=\nu^2\frac{\partial ^2N(x,t)}{\partial x^2}+ \xi^2 N(x,t)+\varphi(x,t), 0\leq\alpha \leq 1
\end{equation}            
with the initial conditions
\begin{equation}
N(x,0)=\delta(x), x\in \Re, N_t(x,0)= 0, \lim_{x\rightarrow \pm \infty}N(x,t) = 0, t>0,
\end{equation}
$\phi(x,t)$ is a nonlinear function belonging to the area of reaction-diffusion. Then for the fundamental solution of (58) subject to the initial conditions (59), there holds the formula 
\begin{eqnarray}
N(x,t)&=&\frac{1}{2\pi\sqrt{(a^2-4b)}}[\int^{+\infty}_{-\infty}exp(-ikx)\nonumber\\
&&\left\{(\lambda+a)E_\alpha(\lambda t^\alpha)-(\mu+a)E_\alpha(\mu t^\alpha)\right\} dk\nonumber\\
&+&\frac{1}{2\pi}\int^t_0\xi^{\alpha-1}\int_{-\infty}^{+\infty}exp(-ikx)\varphi^*(k,t-\xi)\nonumber\\
&&\times [E_{\alpha,\alpha}(\lambda\xi^\alpha)-E_{\alpha,\alpha}(\mu\xi^\alpha)]dkd\xi,
\end{eqnarray}
where $\lambda$ and $\mu$ are the real and distinct roots of the quadratic equation
\begin{equation}
y^2+ay+b=0,
\end{equation} 
given by    
\begin{equation}
\lambda=\frac{1}{2}(-a+\sqrt{(a^2-4b)})\mbox{and}\;\; \mu=\frac{1}{2}(-a-\sqrt{(a^2-4b)}),
\end{equation} 
where $b^2=\nu^2k^2-\xi^2$ and provided that the integral appearing in (60) is convergent.

Next, if we set $\varphi(x,t)=b=0, \eta = 2$, replace $\alpha$ by $2\alpha$,  and $\beta$   by $\alpha$ in  
(58), we then obtain the following interesting result, which includes many known results on fractional telegraph equations including the one recently given by Orsingher and Beghin (2004).\\
\noindent
{\bf Corollary 1.5.} Consider the following reaction-diffusion system
with the initial conditions
\begin{equation}
\frac{\partial ^{2\alpha}N(x,t)}{\partial t^{2\alpha}}+a\frac{\partial^\alpha N(x,t)}{\partial t^\alpha}=\nu^2\frac{\partial^2 N(x,t)}{\partial x^2}+\xi^2 N(x,t), 0 \leq \alpha \leq 1
\end{equation}
with the initial conditions
\begin{equation}
N(x,0)=\delta(x), x\in \Re, N_t (x,0)=0\;\;\lim_{x\rightarrow \pm \infty} N(x,t)=0,t>0,
\end{equation}
Then for the  fundamental solution of (63) subject to the initial conditions (64), there holds the formula 
\begin{eqnarray}
N(x,t)&=&\frac{1}{2\pi\sqrt{(a^2-4b)}}[\int^{+\infty}_{-\infty} exp (-ikx)\nonumber\\
&&\left\{(\lambda+a)E_\alpha(\lambda t^\alpha)-(\mu+a)E_\alpha(\mu t^\alpha)\right\}dk,
\end{eqnarray}
where $\lambda$   and $\mu$  are defined in (62), $b=\nu^2k^2-\xi^2$ and $E_\alpha(x)$ is the Mittag-Leffler function defined by (3). 
      
If we set $\xi^2=0$, then the Corollary 1.5 reduces to the  result, which states that
the reaction-diffusion system (66)
\begin{equation}
\frac{\partial^{2\alpha} N(x,t)}{\partial t^{2\alpha}}+ a\frac{\partial^\alpha N(x,t)}{\partial t^\alpha}=\nu^2 \frac{\partial^2N(x,t)}{\partial x^2},\;\;0\leq \alpha \leq 1
\end{equation}
with the initial conditions
\begin{equation}
N(x,t)=\delta(x), x \in \Re, N_t(x,0)=0,\;\;lim_{x\rightarrow\pm\infty}N(x,t)=0, t>0,
\end{equation}
has the fundamental solution, given by  
\begin{eqnarray}
N(x,t)&=&\frac{1}{2\pi\sqrt{(a^2-4b)}}\int^{+\infty}_{-\infty}exp(-ikx)\nonumber\\
&&\left\{(\lambda+a)E_\alpha(\lambda t^\alpha)-(\mu+a)E_\alpha(\mu t^\alpha)\right\}dk,
\end{eqnarray}
where $\lambda$  and $\mu$ are defined in  (62), $b=\nu^2k^2$ and $E_\alpha(x)$ is the Mittag-Leffler function defined by (3).

The result (68) can be rewritten in the form 
\begin{eqnarray}
&&N(x,t)= \frac{1}{4\pi}\int^{+\infty}_{-\infty}exp(-ikx)\nonumber\\
&&\left\{(1+\frac{a}{\sqrt{(a^2-4\nu^2k^2)}})E_\alpha(\lambda t^\alpha)+\right.\nonumber\\
&&\left.(1-\frac{a}{\sqrt{(a^2-4\nu^2k^2)}})E_\alpha(\mu t^\alpha)\right\}dk,
\end{eqnarray}
where $\lambda$ and $\mu$ are defined in (62) and is the Mittag-Leffler function, defined by (3).

The eq. (69) represents the solution of time-fractional telegraph eq. (66) subject to the initial conditions in (67), recently solved by Orsingher and Beghin (2004). It may be remarked here that the solution as given by Orsingher and Beghin (2004) is in terms of the Fourier transform of the solution in the form given below. It is interesting to observe that the Fourier transform of the solution of the eqs. (66) and (67) can be expressed in the form
\begin{equation}
N^*(x,t)=\frac{1}{2}\left\{(1+\frac{a}{\sqrt{(a^2-4\nu^2k^2)}})E_\alpha(\lambda t^\alpha)+(1-\frac{a}{\sqrt{(a^2-4\nu^2k^2)}})E_\alpha(\mu t^\alpha)\right\},
\end{equation}  
where $\lambda$ and $\mu$ are defined in (62) and $E_\alpha(x)$ is the Mittag-Leffler function defined by (3).
     
Following a similar procedure, the following general theorem can be established which makes use of the general result (24).

\section{General Theorem 2} 

Consider the unified fractional reaction-diffusion equation associated with Caputo derivatives
\begin{eqnarray}
&&a_1\;_0D_t^{\alpha_1} N(x,t)+a_2\;_0D_t^{\alpha_2} N(x,t)+\ldots +a_{n+1}\;_0D_t^{\alpha_{n+1}}N(x,t)\nonumber\\
&=& \nu^2\;_{-\infty}D_x^\mu  N(x,t)+\xi^2 N(x,t)+\varphi(x,t)\\
&&0\leq \alpha_j \leq 1, (j=1,\ldots, n+1),\; \mu>0\nonumber
\end{eqnarray}                                   
with initial conditions 
\begin{equation}
N(x,0)=f(x),\;\;x\in \Re, \;lim_{x\rightarrow \pm \infty} N(x,t)=0,\;\; t>0
\end{equation}
where $\nu^2$ is a diffusion constant, $\xi$ is a constant which  describes the nonlinearity in the system, and $\phi(x,t)$ is nonlinear function which belongs to the area of reaction- diffusion, then there holds the following formula for the solution of (71).
\begin{eqnarray}
N(x,t)&=& \frac{1}{2\pi a_1}\sum^\infty_{m=0}(-1)^m \int^\infty_{-\infty}f^*(k)exp(-ikx)\sum_{r_1+\ldots+r_n=m}\nonumber\\
&=& \frac{m!}{r_1!\ldots r_n!}\left\{\Pi^n_{j=1}\left(\frac{a_{j+1}}{a_1}\right)^{r_j}\right\}\nonumber\\
&\times& t^A\sum^{n+1}_{\omega=1}a_\omega t^{-\alpha_\omega} E^{m+1}_{\alpha_1,A-\alpha_\omega+1}(-\frac{a_0}{a_1}t^{\alpha_1})dk\nonumber\\
&+& \frac{1}{2\pi a_1}\sum^\infty_{m=0}(-1)^m\int^t_0\int^\infty_{-\infty}exp(-ikx)\xi^{A-1}\varphi^*(k,t-\xi)\nonumber\\
&\times&\sum_{r_1+\ldots +r_n=m}\frac{m!}{r_1!\ldots r_n!}\left\{\Pi^n_{j=1}\left(\frac{a_{j+1}}{a_1}\right)^{r_j}\right\}\nonumber\\
&\times& E^{m+1}_{\alpha_1,A}(-\frac{a_0}{a_1} t^{\alpha_1})dkd\xi,
\end{eqnarray}
where $a_0=\nu^2|k|^\mu-\xi^2, Re(\alpha_j)>0, )j=1,\ldots,n); Re(\alpha_{j+1}-\alpha_2)>0$\\ 
for $j=2,\ldots, n, lim_{x\rightarrow\pm \infty} N(x,t)=0; A=\sum^n_{j=1}(\alpha_2-\alpha_{j+1})r_j+\alpha_1(1+m)-\alpha_2m$\\
and provided that the series and integrals in (73) are convergent.

\section{Consequences of General Theorem 2}

When  $\alpha_j =\mu=1/2 (j=1,2,\ldots , n=1)$, the theorem 2 yields \\
{\bf Corollary 2.1.} Consider the unified fractional reaction-diffusion equation associated with Caputo derivatives
\begin{eqnarray} 
&&a_1\;_0D_t^{1/2} N(x,t)+a_2\;_0D_t^{1/2}N(x,t)+\ldots+a_{n+1} \;_0D_t^{1/2} N(x,t)\nonumber\\
&=& \nu^2 \;_{-\infty}D^{1/2}_x N(x,t)+\xi^2 N(x,t)+\varphi(x,t)\\
&& 0\leq \alpha_j \leq 1, (j=1,\ldots, n+1)\nonumber
\end{eqnarray}                        
with initial conditions 
\begin{equation}
N(x, 0) = f(x), (x\in \Re), lim_{x\rightarrow\pm\infty}N(x,t)=0, t>0
\end{equation}
where $\nu^2$  is a diffusion constant, $\xi$  is a constant which  describes the nonlinearity in the system, 
and $\phi(x,t)$ is nonlinear function which belongs  to the area of reaction-diffusion, then there holds the following formula for the solution of (74) and (75).
\begin{eqnarray}
&&N(x,t)=\frac{1}{2\pi a_1}\sum^\infty_{m=0}(-1)^m\int^\infty_{-\infty} f^*(k)exp(-ikx)\sum_{r_1+\ldots+r_n=m}\frac{m!}{r_1!\ldots r_n!}\nonumber\\
&&\left\{\Pi^n_{j=1}\left(\frac{a_{j+1}}{a_1}\right)^{r_j}\right\}
(a_1+\ldots+a_{n+1})E^{m+1}_{1/2,1}(-\frac{bt^{1/2}}{a_1})dk\nonumber\\
&&\sum^\infty_{m=0}(-1)^m\int^t_0 \int^\infty_{-\infty}exp(-ikx)\xi^{-1/2}\varphi^*(k,t-\xi)\nonumber\\
&&\sum_{r_1+\ldots +r_n=m}\frac{m!}{r_1!\ldots r_n!}\left\{\Pi^n_{j=1}\left(\frac{a_{j+1}}{a_1}\right)^{r_j}\right\}+\nonumber\\
&&+\frac{1}{2\pi a_1}E^{m+1}_{1/2,1/2}(-\frac{b}{a_1}t^{1/2})dkd\xi,
\end{eqnarray}               
where $b=\nu^2|k|^{1/2}-\xi^2, a_1\neq 0$, and provided that the series and integrals in eq. (76) are convergent.\\
{\bf Corollary 2.2.} Consider the fractional reaction-diffusion equation 
\begin{eqnarray}
&&a_1\;_0D_t^{\alpha_1}N(x,t)+a_2\;_0D_t^{\alpha_2}N(x,t)+\ldots+a_{n+1}\;_0D_t^{\alpha_{n+1}}N(x,t)\nonumber\\
&=&\nu^2\;_{-\infty}D_x^\mu N(x,t)+\xi^2\;N(x,t)+\varphi(x,t)\\
&& 0\leq a_j\leq 1,(j=1,\ldots, n+1), \mu>0\nonumber
\end{eqnarray}
with initial conditions 
\begin{equation}
N(x,0)=\delta(x), \;(x\in \Re), lim_{x\rightarrow \pm\infty}N(x,t) = 0, t>0,
\end{equation}                                                     
where $\nu^2$ is a diffusion constant, $\xi$ is a constant which describes the nonlinearity in the system, and $\varphi(x,t)$ is nonlinear function which belongs  to the area of reaction-diffusion, then there holds the following formula for the fundamental solution of eqs. (77) and (78).
\begin{eqnarray}
N(x,t) &=& \frac{1}{2\pi\alpha_1}\sum^\infty_{m=0}(-1)^m \int^\infty_{-\infty}exp(-ikx)\sum_{r_1+r_n=m}\frac{m!}{r_1!\ldots r_!}\nonumber\\
&&\left\{\Pi^n_{j=1}\left(\frac{a_{j+1}}{a_1}\right)^{r_j}\right\}
t^A\sum^{n+1}_{\omega=1}t^{-\alpha_\omega}E^{m+2}_{\alpha_1, A-\alpha_\omega+1}(-\frac{a_0}{a_1}t^{\alpha_1})dk\nonumber\\
&+&\frac{1}{2\pi a_1}\sum^\infty_{m=0}(-1)^m \int^t_0\int^\infty_{-\infty}\xi^{A-1} exp(-ikx)\varphi^*(k,t-\xi)\nonumber\\
&&\sum_{r_1+\ldots +r_n=m}\frac{m!}{r_1!\ldots r_n!}\left\{\Pi^n_{j=1}\left(\frac{a_{j+1}}{a_1}\right)^{r_j}\right\}\nonumber\\
&\times& E^{m+1}_{\alpha_1, A}(-\frac{a_0}{a_1}t^{\alpha_1})dkd\xi,
\end{eqnarray}
where $A=\sum^n_{j=1}(\alpha_{j+1}-\alpha_2)r_j+\alpha_1+(\alpha_1-\alpha_2)m, Re(\alpha_j)>0, (j=1,\ldots, n);\\ 
Re (\alpha_{j+1}-\alpha_2)>0$\\

for $j=2,\ldots n, a_0=\nu^2|k|^\mu-\xi^2,\; lim_{x\rightarrow\pm\infty} N(x,t)=0, \;a_1\neq 0$, and provided that the integrals and series in (79) are convergent.\\
On the other hand, if we set $\xi=0$, we arrive at\\ 
{\bf Corollary 2.3.} Consider the fractional reaction-diffusion equation \\
\begin{eqnarray}                                   
&&a_1\;_0D_t^{\alpha_1}N(x,t)+a_2\;_0D_t^{\alpha_2} N(x,t)+\ldots+a_{n+1}\;_0D_t^{\alpha_{n+1}}N(x,t)\nonumber\\
&=& \nu^2\;_{-\infty}D_x^\mu N(x,t)+\varphi(x,t)\\
&&0\leq\alpha\leq 1,0\leq \beta\leq 1,0\leq \gamma\leq 1,\mu>0\nonumber
\end{eqnarray}                  
with initial conditions 
\begin{equation}
N(x,0)=\delta(x), (x\in \Re), lim_{x\rightarrow \pm \infty}N(x,t)=0,t>0,
\end{equation}
where $\nu^2$ is a diffusion constant, and $\phi(x,t)$ is nonlinear function which belongs to the area of reaction-diffusion, then there holds the following formula for the fundamental  solution of eqs. (80) and (81).
\begin{eqnarray}
&&N(x,t)=\frac{1}{2\pi a_1}\sum^\infty_{m=0}(-1)^m \sum_{r_1+\ldots+r_n=m}\frac{m!}{r!_1\ldots r_n!}\left\{\Pi^n_{j=1}\left(\frac{a_{j+1}}{a_1}\right)^{r_j}\right\}\nonumber\\
&&\sum^{n+1}_{\omega=1}a_\omega\int_0^\alpha G_1(x-\tau,t)f(\tau)d\tau\nonumber\\
&&+\frac{1}{2\pi a_1}\sum^\infty_{m=0}(-1)^m\int_0^t(t-\xi)^{A-1}\sum_{r_1+\ldots+r_n=m}\frac{m!}{r_1!\ldots r_n!}\left\{\Pi^n_{j=1}\left(\frac{a_{j+1}}{a_1}\right)^{r_j}\right\}\nonumber\\
&\times&\int^
\alpha_0 G_2(x-\tau,t-\xi)\varphi(\tau,\xi)d\tau d\xi,     
\end{eqnarray}
where, $A=\sum^n_{j=1}(\alpha_{j+1}-\alpha_2)r_j+\alpha_1+(\alpha_1-\alpha_2)m, Re(\alpha_j)>0,\; (j=1, \ldots, n);$\\
$ Re(\alpha_{j+1}-\alpha_2)>0$\\
for $j=2,\ldots, n, a_0=\nu^2|k|^\mu-\xi^2,\;lim_{x\rightarrow \pm \infty}N(x,t)=0,\;\mbox{and}\;\; a_1\neq 0;$
\begin{eqnarray}
G_1(x,t)&=& \frac{t^{A-\alpha_\omega}}{2\pi}\int^\infty_{-\infty}exp(-ikx)E^{m+1}_{\alpha_1,A-\alpha_\omega+1}\left(\frac{\nu^2|k|^\mu t^\alpha}{a_1}\right)dk\nonumber\\
&=&\frac{t^{A-\alpha_\omega}}{\pi\mu(m!)}\int^\infty_0 cos(kx)H^{1,1}_{1,2}\left[\frac{\nu^{2/k}|k|t^{\alpha_1/\mu}}{a_1^{1/\mu}}\left|^{(-m, 1/\mu}_{(0,1/\mu), (\alpha_\omega-\frac{A_1}{\mu},\frac{\alpha_1}{\mu})}\right.\right],\nonumber\\
&=& \frac{t^{A-\alpha_\omega}}{\mu(m!)|x|}\;H^{2,1}_{3,3}\left[\frac{|x|a^{1/\mu}}{\nu^2/\mu_t\alpha_1/\mu}\left|^{(1,1/\mu),(1+A-\alpha_\omega, \alpha_1/\mu), (1,1/2)}_{(1,1),(1+m,1/\mu),(1,1/2)}\right.\right],
\end{eqnarray}
and 
\begin{eqnarray}
G_2(x,t)&=&\frac{t^A}{2\pi}\int^\infty_{-\infty} exp(-ikx) E^{m+1}_{\alpha_1,A+1}\left(-\frac{\nu^2|k|^\mu t^{\alpha_1}}{a_1}\right)dk\nonumber\\
&=& \frac{t^A}{\mu(m)!}\int^\infty_{-\infty} cos(kx)H^{1,1}_{1,2}\left[\frac{\nu^{2/\mu}t^{\alpha_1/\mu}|k|}{a_1^{1/\mu}}\left|^{(-m,1/\mu)}_{(0,1/\mu), (1-A,\alpha_1/\mu)}\right.\right]\nonumber\\
&=& \frac{t^A}{\mu(m)!|x|}H^{2,1}_{3,3}\left[\frac{a_1^{1/\mu}|x|}{\nu^2/\mu}\left|^{(1,1/\mu),(A,\alpha_1/\mu),(1,1/2)}_{(1,1),(1+m,1/\mu), (1,1/2)}\right.\right],
\end{eqnarray}     
 
provided that the series and integral in eq. (82) are convergent.\par
\bigskip
\noindent
{\bf Remark.} It is interesting to observe that the method employed for deriving the solution of the equations (41) and (42) as well as (71) and (72) in the space $=LF=L(R_+)\times F(R)$ can also be applied  in the space $LF'=L'(R_+)\times F'$, where $F'=F'(r)$  is the space of Fourier transform of  generalized functions .As an illustration, we can choose  $F' = S'$ or $F'=D'$. The Fourier transforms in $S'$ and $D'$ are introduced by Schwartz and Gelfand and Shilov,  respectively. $S'$ is the dual of the  space $S$, which  is the space of all infinitely differentiable functions which, together with their derivatives approach zero more rapidly than any power of  $1/|x|$ as $|x|\rightarrow \infty$ (Gelfand and Shilov, 1964, p. 16). Here $D'$ is the dual of the space $D$  which consists of all smooth functions with compact supports (Brychkov and Prudnikov, 1989, p. 3). For further details, the reader is referred to the monographs written by Gelfand and Shilov (1964) and Brychkov and Prudnikov et al. (1989), if we replace  the Laplace and Fourier transforms in eqs. (10) and (12) by the corresponding Laplace and Fourier transform of the generalized functions.\par
\bigskip
\noindent
{\bf Acknowledgment} The authors would like to thank the Department of Science and Technology, Government of India, New Delhi, for the financial assistance for this work under project No. SR/S4/MS:287/05 which enabled this collaboration possible.\par
\bigskip
\noindent
{\bf  References}\\
Abramowitz, M. and Stegun, I.A.: 1968, {\it Handbook of Mathematical Functions with Formulas, Graphs, and Mathematical Tables}, Applied Math. Series 55, 7th Printing, National Bureau of Standards, Washington, DC.\par
\smallskip
\noindent
Brychkov, Yu.A. and Prudnikov, A.P.: 1989, {\it Integral Transforms of Generalized Functions}, Gordon and Breach, New York.\par
\smallskip
\noindent
Caputo, M.: 1969, {\it Elasticita e Dissipazione}, Zanichelli, Bologna.\par
\smallskip
\noindent
Compte, A.: 1996, Stochastic foundations of fractional dynamics, {\it Physical Review E}, {\bf 53}, 4191-4193.\par
\smallskip
\noindent
Doetsch, G.: 1956, {\it Anleitung zum Praktischen Gebrauch der Laplace-\\
Transformation}, Oldenbourg, Munich.\par
\smallskip
\noindent
Dzherbashyan, M.M.: 1966, {\it Integral Transforms and Representation of Functions in Complex Domain} (in Russian), Nauka, Moscow.\par
\smallskip
\noindent 
Dzherbashyan, M.M.: 1993, {\it Harmonic Analysis and Boundary Value Problems in the Complex Domain}, Birkhaeuser-Verlag, Basel and London.\par
\smallskip
\noindent
Erd\'{e}lyi, A., Magnus, W., Oberhettinger, F., and Tricomi, F.G.: 1953, {\it Higher Transcendental Functions}, Vol. {\bf 1}, McGraw-Hill, New York, Toronto.\par
\smallskip
\noindent
Erd\'{e}lyi, A., Magnus, W., Oberhettinger, F., and Tricomi, F.G.: 1954, {\it Tables of Integral Transforms}, Vol. {\bf 2}, McGraw-Hill, New York, Toronto.\par
\smallskip
\noindent
Erd\'{e}lyi, A., Magnus, W., Oberhettinger, F., and Tricomi, F.G.: 1955, {\it Higher Transcendental Functions}, Vol. {\bf 3}, McGraw-Hill, New York, Toronto.\par
\smallskip
\noindent
Gelfand, I.M. and Shilov, G.F.: 1964, {\it Generalized Functions}, Vol. {\bf 1}, Academic Press, London.\par
\smallskip
\noindent
Kilbas, A.A. and Saigo, M.:  2004, {\it H-transforms: Theory and Applications}, CRC Press, New York.\par
\smallskip
\noindent
Kilbas, A.A., Saigo, M., and Saxena, R.K.: 2002, Solution of Volterra integro-differential equations with generalized Mittag-Leffler function in the kernels, {\it Journal of Integral Equations}, {\bf 14}, 377-396.\par
\smallskip
\noindent 
Kilbas, A.A., Saigo, M., and Saxena, R.K.: 2004, Generalized Mittag-Leffler function and generalized fractional calculus, {\it Integral Transforms and Special Functions}, {\bf 15}, 31-49.\par
\smallskip
\noindent
Manne, K.K., Hurd, A.J., and Kenkre, V.M.: 2000, Nonlinear waves in reaction-diffusion systems: The effect of transport memory, {\it Physical Review E}, 
{\bf 61}, 4177-4184.\par
\smallskip
\noindent
Metzler, R. and Klafter, J.: 2000, The random walk's guide to anomalous diffusion: a fractional dynamics approach, {\it Physics Reports}, {\bf 339}, 1-77.\par
\smallskip
\noindent
Miller, K.S. and Ross, B.: 1993, {\it An Introduction to the Fractional Calculus and Fractional Differential Equations}, John Wiley and Sons, New York.\par
\smallskip
\noindent
Oldham, K.B. and Spanier, J.: 1974, {\it The Fractional Calculus: Theory and Applications of Differentiation and Integration of Arbitrary Order}, Academic Press, New York.\par
\smallskip
\noindent
Orsingher, E. and Beghin, L.: 2004, Time-fractional telegraph equations and telegraph processes with Brownian time, {\it Probability Theory and Related Fields} 
{\bf 128}, 141-160.\par
\smallskip
\noindent
Podlubny, I.: 1999, {\it Fractional Differential Equations}, Academic Press, San Diego.\par
\smallskip
\noindent
Prabhakar, T.R.: 1971, A singular integral equation with generalized Mittag-Leffler function in the kernel, {\it Yokohama Mathematical Journal}, {\bf 19}, 7-15.\par
\smallskip
\noindent
Prudnikov, A.P., Brychkov, Yu.A., and Marichev, O.I.: 1989, {\it Integrals and Series, More Special Functions}, Vol. {\bf 3}, Gordon and Breach, New York.\par
\smallskip
\noindent
Samko, S.G., Kilbas, A.A., and Marichev, O.I.: 1990.: {\it Fractional Integrals and Derivatives: Theory and Applications}, Gordon and Breach, New York.\par
\smallskip
\noindent 
Saxena, R.K., Mathai, A.M., and Haubold, H.J.: 2002, On fractional kinetic equations, {\it Astrophysics and Space Science}, {\bf 282}, 281-287.\par
\smallskip
\noindent 
Saxena, R.K., Mathai, A.M., and Haubold, H.J.: 2004, On  generalized fractional kinetic equations, {\it Physica A}, {\bf 344}, 657-664.\par
\smallskip
\noindent
Saxena, R.K., Mathai, A.M., and Haubold, H.J.: 2004a, Unified fractional kinetic equation and a fractional diffusion equation, {\it Astrophysics and Space Science}, {\bf 290}, 299-310\par
\smallskip
\noindent
Saxena, R.K., Mathai, A.M., and Haubold, H.J.: 2004b, Astrophysical thermonuclear functions for Boltzmann-Gibbs statistics and Tsallis statistics, {\it Physica A}, {\bf 344}, 649-656.\par
\smallskip
\noindent
Saxena, R.K., Mathai, A.M., and Haubold, H.J.: 2006a, Fractional reaction-diffusion equations, {\it Astrophysics and Space Science}, {\bf 305}, 289-296.\par
\smallskip
\noindent
Saxena, R.K., Mathai, A.M., and Haubold, H.J.: 2006b, Reaction-diffusion systems and nonlinear waves, {\it Astrophysics and Space Science}, {\bf 305}, 297-303.\par
\smallskip
\noindent
Saxena, R.K., Mathai, A.M., and Haubold, H.J.: 2006c, Solution of generalized fractional reaction-diffusion equations, 
{\it Astrophysics and Space Science}, {\bf 305}, 305-313.\par
\smallskip
\noindent
Strier, D.E., Zanette, D.H., and Wio, H.S.: 1995, Wave fronts in a bistable reaction-diffusion system with density-dependent diffusivity, 
{\it Physica A}, {\bf 226}, 310-323.\par
\smallskip
\noindent
Wiman, A.: 1905, Ueber den Fundamentalsatz in der Theorie der Funktionen E(x), {\it Acta Mathematica}, {\bf 29}, 191-201.
\end{document}